\documentclass[12pt]{article}
\usepackage{amsthm}
\usepackage{amsmath,amsfonts,amssymb,epsfig,amsbsy}
\begin{document}
\def\Vect{\operatorname{Vect}}
\def\Pol{\operatorname{Pol}}
\def\Hom{\operatorname{Hom}}
\def\codim{\operatorname{codim}}
\def\tame{{\operatorname{tame}}}
\def\L{\operatorname{L}}
\def\S{\operatorname{S}}
\def\C{\operatorname{C}}
\def\gp{\operatorname{gp}}
\def\GL{\operatorname{GL}}
\def\Col{\operatorname{Col}}
\def\join{\operatorname{join}}
\def\N{\operatorname{N}}
\def\Im{\operatorname{Im}}
\def\conv{\operatorname{conv}}
\def\htt{\operatorname{ht}}
\def\ker{\operatorname{ker}}
\newtheorem{conjecture}{Conjecture}
\newtheorem{proposition}{Proposition}
\newtheorem{theorem}{Theorem}
\newtheorem{example}{Example}
\newtheorem{lemma}{Lemma}
\newtheorem{remark}{Remark}

\title{TAME HOMOMORPHISMS OF POLYTOPAL RINGS}
\author{Viveka Erlandsson\\
Department of Mathematics\\
City University of New York, Graduate Center\\
365 Fifth Avenue\\
New York, NY 10016\\
email: verlandsson@gc.cuny.edu}
\date{}
\maketitle

\begin{abstract}
The object of this paper is the tameness conjecture introduced
in~\cite{BGA}, which describes an arbitrary graded $k$-algebra
homomorphism of polytopal rings. We will give further evidence of
this conjecture by showing supporting results concerning joins,
multiples and products of polytopes.
\end{abstract}



\section{Introduction}
This paper studies the category of polytopal algebras over a field
$k$, denoted $\Pol(k)$. We investigate the concept of
\emph{tameness}, as introduced in~\cite{BGA}, where it is
conjectured that every graded $k$-algebra homomorphism is tame. In
short, such a homomorphism is called tame if it can be obtained by a
composition of some standard homomorphisms, defined in~\cite{BGA}.

In this work we will show the following: If graded homomorphisms
between two polytopal rings are tame, then graded homomorphisms from
rings obtained by taking multiples, joins, and products of the
underlying polytopes are also tame. Thus, we extend the class of
polytopal rings on which graded homomorphisms are tame, giving
further evidence in support of the mentioned conjecture.

The objects of $\Pol(k)$, the \emph{polytopal (monoid) rings}, are
defined as follows. Let $P$ be a convex lattice polytope in
$\mathbb{R}^{n}$. Let $\L(P)$ denote the lattice points in $P$, i.e.
$\L(P)=P\cap\mathbb{Z}^{n}$, and let $\S(P)$ be the additive monoid
of $\mathbb{Z}^{n+1}$ generated by $\{(x,1)\,|\,x\in \L(P)\}$. A
lattice point in $\S(P)$ can be represented as a monomial in $n+1$
variables by identifying its coordinate vector with the monomial's
exponent vector. The degree of a monomial is in this case the last
component of its exponent vector. The polytopal ring $k[P]$ is the
monoid ring of $\S(P)$ with coefficients in $k$. $k[P]$ is a graded
ring generated by its degree $1$ monomials. These monomials
correspond bijectively to the lattice points in $P$. The generators'
relations are the binomial relations representing the affine
dependencies in $\L(P)$.

A homomorphism in $\Pol(k)$ is a homomorphism of two polytopal rings
as $k$-algebras preserving the grading. Such a homomorphism also has
to preserve the binomial relationships among the generators.
Therefore, $\Hom(k[P], k[Q])$, which is the set of all graded
homomorphism between the two polytopal rings, is the zero set of a
system of polynomials. Hence $\Hom(k[P], k[Q])$ gives rise to a
Zariski closed set in the space of matrices $M_{mn}(k)$, where
$m=\#L(Q)$ and $n=\#L(P)$, called the \emph{Hom-variety}. We are
aiming to describe the structure of the Hom-variety. The tameness
conjecture is a geometric description of this set.

A tame homomorphism is a graded k-algebra homomorphism which can be
obtained by a composition of four standard homomorphisms: polytope
changes, homothetic blow-ups, Minkowski sums, and free extensions.

For the following definitions, assume $f: k[P]\rightarrow k[Q]$ is a
graded homomorphism and $P,Q\subset \mathbb{R}^{n}$ are lattice
polytopes.

First, suppose $P'$ and $Q'$ are lattice polytopes such that
$P'\subset P$, $Q'\supset Q$, and $f(k[P'])\subset k[Q']$. This
gives rise to a new homomorphism $f': k[P']\rightarrow k[Q']$
obtained from $f$ in a natural way. Also, two polytopes $\tilde{P}$
and $\tilde{Q}$ that are isomorphic to $P$ respectively $Q$ as
lattice polytopes, result in the homomorphism $\tilde{f}:
k[\tilde{P}]\rightarrow k[\tilde{Q}]$, induced by $f$. The
homomorphisms obtained these ways are called \emph{polytope
changes}.

Second, consider the normalization of $S(P)$, defined by
\[
\overline{S(P)}=\{x\in \gp(S(P))\thinspace | \thinspace x^{m}\in
S(P) \mbox{ for some } m\in \mathbb{N}\}\,
\]
and the normalization of $S(Q)$, defined similarly. The set
$\gp(S(P))$ denotes the group of differences of $S(P)$. $S(P)$ is
normal if and only if $S(P)=\overline{S(P)}$. It is well known that
$k[\overline{P}]:=k[\overline{S(P)}]=\overline{k[P]}$. If there are
no monomials in the kernel of $f$, i.e. $\ker(f)\cap
S(P)=\emptyset$, then $f$ extends uniquely to the homomorphism
$\overline{f}:k[\overline{P}]\rightarrow k[\overline{Q}]$ defined by
\[
\overline{f}(x)=\frac{f(y)}{f(z)} \quad\mbox{where}\medspace\,\,
x=\frac{y}{z},\medspace x\in
\overline{S(P)},\medspace\mbox{and}\medspace\, y,z\in S(P).
\]
In fact, $f(y)/f(z)$ belongs to $k[\overline{Q}]$ since
$x\in\overline{S(P)}$ implies $x^{c}\in S(P)$, for some natural
number $c$, which in turn implies $f(x^{c})=(f(y)/f(z))^{c}\in
k[Q]$. Let $k[\overline{P}]_{c}$ be the subring of $k[\overline{P}]$
generated by the homogeneous components of degree $c$ in
$\overline{S(P)}$ (and similarly for $k[\overline{Q}]_{c}$). Note
that $k[\overline{P}]_{c}\simeq k[cP]$ and
$k[\overline{Q}]_{c}\simeq k[cQ]$ in a natural way. Since $f$ is
graded, $\overline{f}$ restricts to the graded homomorphism
$f^{(c)}: k[cP]\rightarrow k[cQ]$, which we call the
\emph{homothetic blow-up} of $f$.

Third, consider having two graded homomorphisms $f,g:k[P]\rightarrow
k[Q]$. Let $\N(f(x))$ denote the Newton polytope of $f(x)$, i.e. the
convex hull of the support monomials in $f(x)$ (and similarly for
$\N(g(x))$). Assume $f$ and $g$ satisfy
\[
\N(f(x))+\N(g(x))\subset Q  \,\, \mbox{for all}\, x\in L(P)
\]
where $+$ denotes the Minkowski sum in $\mathbb{R}^{n}$. We have
$z^{-1}f(x)g(x)\in k[Q]$ where $z=(0,0,...,1)$, and thus
\[
f\star g: k[P]\rightarrow k[Q] \, \mbox{ such that }\, f\star
g(x)=z^{-1}f(x)g(x)
\]
for all $x\in L(P)$ defines a new graded homomorphism, called the
\emph{Minkowski sum} of $f$ and $g$.

Lastly, assume $P$ is a pyramid with base $P_{0}$ and vertex $v$
such that $L(P)=\{v\}\cup L(P_{0})$ and that we are given a
homomorphism $f_{0}: k[P_{0}]\rightarrow k[Q]$. This means $k[P]$ is
a polynomial extension of $k[P_{0}]$. Thus, $f_{0}$ extends to a
homomorphism $f:k[P]\rightarrow k[Q]$ by letting $f(v)=q$ for any
$q\in k[Q]$ and $f(x)=f_{0}(x)$ for all $x\in \L(P_{0})$. The
homomorphism $f$ is called a \emph{free extension} of $f_{0}$.

The tameness conjecture, as found in~\cite{BGA}, states the
following.

\begin{conjecture}\emph{(W. Bruns, J. Gubeladze)}
Every homomorphism in $\Pol(k)$ is obtained by a sequence of taking
free extensions, Minkowski sums, homothetic blow-ups, polytope
changes and compositions, starting from the identity mapping $k\to
k$. Moreover, there are normal forms of such sequences.
\end{conjecture}

Certain subvarieties of the Hom-variety have already been described.

In~\cite{BGG} the subvariety corresponding to automorphisms has been
described completely, showing that every automorphism is a
composition of some basic automorphisms. This can be viewed as a
polytopal generalization of the linear algebra fact that every
invertible matrix can be written as the product of elementary,
permutation, and diagonal matrices. The result is stronger than the
notion of tameness, thus implying that automorphisms are tame.

In~\cite{BGR} the variety corresponding to codimension-1 retractions
of polygons has been described. It is conjectured that the result
generalizes to arbitrary dimensions of the polytope. This result is
also stronger than tameness; hence, such retractions are tame.

In showing that these results imply tameness of the corresponding
morphisms, some other classes of homomorphisms are shown to be tame
in~\cite{BGA}. Among these classes, the following two will be used.
First, homomorphisms respecting monomial structures are tame.
Second, homomorphisms from $k[c\Delta_{n}]$ (where $\Delta_{n}$ is
the $n$-simplex, $c$, $n\in\mathbb{N}$) are tame. We will also use
the notion of face retractions, which are idempotent endomorphisms
of $k[P]$ defined by
\[
\pi_{F}(x)=\begin{cases}x, &\text{if $x\in F$}\\0, &\mbox{if
$x\notin F$}\end{cases}
\]
for all $x\in\L(P)$, where $F$ is a face of $P$ (note that face
retractions are tame since they respect the monomial structures).

\emph{Acknowledgement:}  This work was suggested by Dr. Joseph
Gubeladze at San Francisco State University.

\section{Main Results}
The results of this work concern joins, multiples, and Segre
products.

The join of two lattice polytopes $P\subset\mathbb{R}^{n}$ and
$Q\subset\mathbb{R}^{m}$ of dimension $n$ respectively $m$ is a
subset of $\mathbb{R}^{n+m+1}$. Consider the embeddings $\iota_{1}$
and $\iota_{2}$ of $\mathbb{R}^{n}$ and $\mathbb{R}^{m}$ into
$\mathbb{R}^{n+m+1}$ defined as
\[
\iota_{1}: (x_{1},\ldots,x_{n})\mapsto
(x_{1},\ldots,x_{n},0,\ldots,0) \qquad \mbox{for all} \,\,
(x_{1},\ldots,x_{n})\in\mathbb{R}^{n}
\]
\[
\,\quad\iota_{2}: (y_{1},\ldots,y_{m})\mapsto
(0,\ldots,0,y_{1},\ldots,y_{m},1) \quad \mbox{for all} \,\,
(y_{1},\ldots,y_{m})\in\mathbb{R}^{m}\,\,\,
\]
The join of $P$ and $Q$ is defined as the convex hull of the image
of $P$ under $\iota_{1}$ and the image of $Q$ under $\iota_{2}$,
i.e. $\join(P,Q)=\conv(\iota_{1}(P),\iota_{2}(Q))$.

We denote the $c^{th}$ multiple of $P$ by $cP$, i.e. $cP=\{cx\, |\,
x\in P\}$. Since the lattice points in $S(P)$ sitting on height $c$
may correspond to a proper subset of the lattice points in $cP$,
$k[cP]$ is in general an overring of the $c^{th}$ Veronese subring
of $k[P]$ (the ring generated by the homogeneous components of
degree $c$ in $S(P)$). The two rings coincide when $P$ is normal.

The Segre product of $k[P]$ and $k[Q]$ is $k[P\times Q]$ where
$P\times Q = \{(x,y) | \, x\in P, y\in Q\}$.

\begin{theorem}
If every graded homomorphism from $k[P]$ respectively $k[Q]$ is
tame, then
\begin{itemize}
\item[(a)] every graded homomorphism from $k[\join(P,Q)]$ is tame,
\item[(b)] every graded homomorphism from $k[cP]$, where $c\in$ $\mathbb{N}$, is
tame,
\item[(c)] every graded homomorphism from $k[P\times Q]$ is tame.
\end{itemize}
\end{theorem}

To prove Theorem $1$, the following lemma will be used. The lemma is
stated and proved in the case of Segre products, but similar
arguments hold true for joins and multiples.

\begin{lemma}
Let $f$ be a graded homomorphism from $k[P\times Q]$ and assume the
hypothesis of Theorem 1. To show that $f$ is tame we can without
loss of generality assume that $\ker(f)\cap S(P\times Q)=\emptyset$.
\end{lemma}
\begin{proof}
Let $f:k[P\times Q]\to k[R]$ be a graded homomorphism (where $R$ is
some lattice polytope). First we observe that the hypothesis of
Theorem 1 descends to the faces of the polytopes. That is, if every
graded homomorphism from $k[P]$ is tame and $F$ is a face of $P$,
then any graded homomorphism from $k[F]$ is tame. Note that such a
homomorphism $f':k[F]\to k[R]$ is the composition of the
homomorphisms
\[
k[F]\xrightarrow{\text{$x\mapsto
x$}}k[P]\xrightarrow{\text{$x\mapsto
\pi_{F}(x)$}}k[F]\xrightarrow{\text{$x\mapsto f'(x)$}}k[R]
\]
The first map is tame since it maps monomials to monomials. The
composition of the last two maps is tame since it is a homomorphism
from $k[P]$, tame by assumption.

Now, assume there exist a monomial $m\in \ker(f)$. Since $\ker(f)$
is a prime ideal, the ideal $I=(\ker(f)\cap S(P\times Q))\subset
k[P\times Q]$ is a monomial prime ideal containing $m$. However,
monomial prime ideals are exactly the kernels of face
retractions~\cite{Book}. Thus, $f$ is a composition of a face
retraction and a map $g:k[F']\to k[R]$ where $F'$ is a face of
$P\times Q$. Face retractions are known to be tame. Hence, since a
composition of tame homomorphisms is tame, tameness of $g$ will
imply tameness of $f$. Note that there are no monomials in the
kernel of $g$. Also, a face of $P\times Q$ is of the form $P'\times
Q'$ where $P'$ is a face of $P$ and $Q'$ a face of $Q$~\cite{Z}.

Thus, Theorem $1(c)$ is proved by establishing tameness of every
graded homomorphism $g$ from $k[P'\times Q']$ such that $\ker(g)\cap
S(P'\times Q')=\emptyset$, given that homomorphisms from $k[P']$ and
$k[Q']$ are tame.
\end{proof}

From here on, it is assumed that there are no monomials in the
kernel of the homomorphisms in question.

\section{Translation into Discrete Objects}\label{S:discrete}
A graded homomorphism in $\Hom(k[P],k[Q])$ can be viewed as an
affine map from $\L(P)$ to $\mathbb{Z}_{+}^{l}$, and vice versa. In
this case, ``affine" refers to a map admitting an affine extension
to the corresponding affine hulls. By viewing a graded homomorphisms
as an affine map, polynomials in a polytopal ring translate into
discrete objects, which will aid in proving Theorem $1$.

Suppose $P$ and $Q$ are lattice polytopes such that $\L(P)\subset
\mathbb{Z}^{d}_{+}$ and $\L(Q)\subset \mathbb{Z}^{e}_{+}$ (this can
always be assumed by a polytope change). Let $f: k[P]\to k[Q]$ be a
graded homomorphism. $\L(Q)$ is a subset of $\{X_{1}^{a_{1}}\cdot
\cdot \cdot X_{e}^{a_{e}} Z | a_{i}\in\mathbb{Z}_{+}\}$ and thus the
polynomials $\varphi_{x}=f(x)Z^{-1}$ belong to $k[X_{1}, ...,
X_{e}]$ for all $x\in L(P)$. Clearly, since $f$ respects the
binomial relations in $L(P)$, so does the map $x\mapsto\varphi_{x}$.

The polynomial ring $k[X_{1}, ..., X_{e}]$ is a unique factorization
domain. Therefore $(k[X_{1}, ..., X_{e}] \setminus \{0\}) /k^{*}$ is
a free commutative monoid. Hence, there is a subset
$\mathcal{P}\subset (k[X_{1}, ..., X_{e}]\setminus \{0\}) /k^{*}$ of
irreducibles such that for each class $[\varphi_{x}]\in (k[X_{1},
..., X_{e}]\setminus \{0\}) /k^{*}$ we have
$[\varphi_{x}]=P_{1}^{a_{1}} \cdot \cdot \cdot P_{l}^{a_{l}}$ for
some class $P_{i}\in \mathcal{P}$, $l\in \mathbb{N}$, and uniquely
determined $a_{i}\geq 0$. Since only finitely many irreducibles are
needed to represent any $[\varphi_{x}]$ there exists an $l\in
\mathbb{N}$ such that $[\varphi_{x}]\in \{P_{1}^{a_{1}}\cdot \cdot
\cdot P_{l}^{a_{l}} | P_{i}\in \mathcal{P}, a_{i}\in
\mathbb{Z}_{+}\}$ for all $x\in \L(P)$. Note that
$\{P_{1}^{a_{1}}\cdot \cdot \cdot P_{l}^{a_{l}} | P_{i}\in
\mathcal{P}, a_{i}\in \mathbb{Z}_{+}\}$ is isomorphic to
$\mathbb{Z}^{l}_{+}$. Each polynomial $f(x)$ with $x\in L(P)$
therefore gives rise to an integral vector by the correspondence
\[
f(x)\mapsto\ [\varphi_{x}]=P_{1}^{a_{1}} \cdot \cdot \cdot
P_{l}^{a_{l}}\mapsto (a_{1}, \ldots, a_{l})
\]
Consequently, $f$ gives rise to an affine map from $\L(P)$ to
$\mathbb{Z}^{l}_{+}$ which respects the binomial relations in
$\L(P)$.

Conversely, suppose $\alpha$ is an affine map $\L(P)\to
\mathbb{Z}^{l}_{+}$. For each $x\in \L(P)$,
$\alpha(x)\in\mathbb{Z}_{+}^{l}$ gives rise to a polynomial $p_{x}$
in $k[X_{1},\ldots,X_{e}]$ under the correspondence
\[
\alpha(x)=(a_{1},\ldots,a_{l})\mapsto P_{1}^{a_{1}}\cdots
P_{l}^{a_{l}}=p_{x}.
\]
Therefore, for each $x\in L(P)$, $\varphi_{x}=t_{x}p_{x}$ for some
nonzero scalar $t_{x}$. Note that the scalars $t_{x}$ are clearly
subject to the same binomial relations as the lattice points in $P$.
Letting $p'_{x}=t_{x}p_{x}$ for all $x\in L(P)$, results in
$f(x)=p'_{x}Z$, thus recovering the homomorphism $f$ from $\alpha$.

Also, $\alpha$ can be extended to an affine integral map $P\to
\mathbb{R}_{+}^{l}$. If $\L(P)=\{x_{0},\ldots,x_{n}\}$ then $P=\conv
(x_{0},\ldots,x_{n})$ and any $x\in P$ can be represented as
$x=c_{0}x_{0}+\cdots+c_{n}x_{n}$ for some nonnegative real numbers
$c_{i}$ satisfying $\sum_{i=0}^{n}{c_{i}}=1$. Hence, if $f(x_{i})$
corresponds to $a_{i}\in\mathbb{Z}^{l}_{+}$ for each
$i\in\{0,1,\ldots,n\}$, $f(x)$ corresponds to
$c_{0}a_{0}+\cdots+c_{n}a_{n}\in\mathbb{R}_{+}^{l}$. The
correspondence is well-defined since $f$ respects the binomial
relations in $\L(P)$.

\section{Joins}\label{S:veronese}
The proof of Theorem $1(a)$ is mainly derived from well-known
properties of joins. It is known that there are no relations between
the lattice points coming from $P$ and the lattice points coming
from $Q$, in $\join(P,Q)$. Any two polytopal rings $k[P]$ and $k[Q]$
therefore satisfy $k[\join(P,Q)]\simeq k[P]\otimes k[Q]$. Also, by
definition, there are no new lattice points in $\join(P,Q)$ since if
$x\in\L(\join(P,Q))$ then $x\in\L(\iota_{1}(P))\simeq\L(P)$ or $x\in
\iota_{2}(\L(Q))\simeq\L(Q)$. Hence, two graded homomorphisms
$f:k[P]\to k[L]$ and $g:k[Q]\to k[L]$ define a new graded
homomorphism $F:k[\join(P,Q)]\to k[L]$ by letting
\[
F(x)=
\begin{cases}
f(x), &\text{if $x \in \L(P)$}\\
g(x) &\text{if $x \in \L(Q) $}
\end{cases}
\quad\mbox{for all}\,\, x\in \L(\join(P,Q))
\]
Conversely every homomorphism $F:k[\join(P,Q)]\to k[L]$ is
necessarily of this form.

Theorem $1(a)$ is proved once it is shown that this ``pasting" of
the two tame homomorphisms is also tame.

\begin{proof}[Proof of Theorem 1(a)]
Suppose $P$ and $Q$ are lattice polytopes. Let $F: k[\join(P,Q)]\to
k[R]$ be a graded homomorphism where $R$ is a lattice polytope in
$\mathbb{R}^{d}$. Assume $Z=(0,\ldots,0,1)\in \L(R)$, after a
polytope change. Further assume $\L(P)=\{x_{0},\ldots,x_{n}\}$ and
$\L(Q)=\{y_{0},\ldots,y_{m}\}$. There exists homomorphisms $f:
k[P]\to k[R]$ and $g: k[Q]\to k[R]$ such that $F(x_{i})=f(x_{i})$
for $i\in \{0,1,\ldots,n]$ and $F(y_{j})=g(y_{j})$ for
$j\in\{0,1,\ldots,m]$.

Consider the homomorphisms
\[
a_{P}: k[\join(P,Q)]\to k[R] \,\,\mbox{such that}\,\,
a_{P}(x_{i})=f(x_{i}),\,\, a_{P}(y_{j})=Z
\]
\[
a_{Q}: k[\join(P,Q)]\to k[R] \,\,\mbox{such that}\,\,
a_{Q}(x_{i})=Z,\,\, a_{Q}(y_{j})=g(y_{j})
\]
for $i\in\{0,1,\ldots,n\}$ and $j\in\{0,1,\ldots,m\}$. Note that
$a_{P}$ factors through the homomorphism $k[\join(P,Q)]\to
k[\join(P,y_{0})]$ such that $x_{i}\mapsto x_{i}$ and $y_{j}\mapsto
y_{0}$ and the homomorphism $k[\join(P,y_{0})]\to k[R]$ such that
$x_{i}\mapsto f(x_{i})$ and $y_{0}\mapsto Z$. The first map is tame
since it maps monomials to monomials and the second is tame since it
is a free extension of $f$ which is tame by assumption. Thus $a_{P}$
is tame. Similarly $a_{Q}$ factors through the homomorphisms
$k[\join(P,Q)]\to k[\join(x_{0},Q)]$ such that $y_{j}\mapsto y_{j}$
and $x_{i}\mapsto x_{0}$ and $k[\join(x_{0},Q)]\to k[R]$ such that
$y_{j}\mapsto f(y_{j})$ and $x_{0}\mapsto Z$. These maps are tame
for the same reasons and hence $a_{Q}$ is also a tame homomorphism.

$F$ is obtained by the Minkowski sum of $a_{P}$ and $a_{Q}$ since
for all $v\in\L(\join(P,Q))$
\[
(a_{P}\star a_{Q})(v)=a_{P}(v)a_{Q}(v)Z^{-1}=
\begin{cases}
f(v), &\text{if $v \in\L(P)$}\\
g(v), &\text{if $v \in\L(Q) $}
\end{cases}
\]
which is precisely $F$. $F$ is tame as desired.
\end{proof}

\section{Multiples}
Graded homomorphisms from $k[c\Delta_{n}]$, where $c\Delta_{n}$ is
the $c^{th}$ multiple of the $n$-simplex, are tame, as shown
in~\cite{BGA}. A similar argument is used to prove its
generalization concerning the $c^{th}$ multiple of a general lattice
polytope, stated as Theorem $1(b)$. The following lemma is used.

\begin{lemma}
Suppose $P$ is a lattice polytope and $\alpha:
\L(cP)\to\mathbb{Z}_{+}^{l}$ is an affine map. Then there exists a
vector $v \in \mathbb{Z}_{+}^{l}$ and an affine map $\beta:
\L(P)\to\mathbb{Z}_{+}^{l}$ such that $\alpha(cx)=v+c\beta(x)$ for
all $x\in\L(P)$.
\end{lemma}

\begin{proof}[Proof of Lemma 2]
Assume $\L(P)=\{x_{0}, ...\, ,x_{n}\}$. Suppose
$\alpha(cx_{i})=(a_{i1},...,a_{il})$ for $i\in\{0,1,\ldots,n\}$. Let
$v=(\min\{a_{i1}\}_{i=0}^{n},...,\min\{a_{il}\}_{i=0}^{n})\in
\mathbb{Z}_{+}^{l}$. We will show that $\alpha(cx_{i})-v$ are
$c^{th}$ multiples of integral vectors for all $i$.

Let $i\in\{0,1,\ldots,n\}$. Note that for all $k\in\{1,2,\ldots,l\}$
the $k^{th}$ component of either $\alpha(cx_{i})-v$ or
$\alpha(cx_{j})-v$, for some $j\neq i$, is $0$. If the $k^{th}$
component of $\alpha(cx_{i})-v$ is $0$ we are done.

Assume the $k^{th}$ component of $\alpha(cx_{j})-v$ is $0$ for some
$j\neq i$. Since $\alpha$ is an affine map,
$\alpha(cx_{i})-\alpha(cx_{j})=c(\alpha(x_{i})-\alpha(x_{j}))$ is a
$c^{th}$ multiple of an integral vector. But
$\alpha(cx_{i})-\alpha(cx_{j})=(\alpha(cx_{i})-v)-(\alpha(cx_{j})-v)$.
Since the $k^{th}$ component of $\alpha(cx_{j})-v$ is $0$, the
$k^{th}$ component of $\alpha(cx_{i})-\alpha(cx_{j})$ is the
$k^{th}$ component of $\alpha(cx_{i})-v$. Consequently, the $k^{th}$
component is a $c^{th}$ multiple of an integer and it is clearly
nonnegative by the choice of $v$.

Hence $\alpha(cx_{i})-v$ are $c^{th}$ multiples of integral vectors
for all $i$. Denote $\frac{1}{c}(\alpha(cx_{i})-v)$ by
$\beta(x_{i})$ for each $i$. This results in the affine map $\beta:
\L(P)\to\mathbb{Z}_{+}^{l}$ such that $\alpha(cx)=v+c\beta(x)$ for
all $x\in\L(P)$, as desired.
\end{proof}

Note that the positivity condition in Lemma 2 is needed so that the
integral vectors can be identified with polynomials through the
process explained in Section 3.

\begin{proof}[Proof of Theorem 1(b)]
Suppose $\L(P)=\{x_{0},...,x_{n}\}$, $\L(cP) \subset
\mathbb{Z}_{+}^{d}$ and $\L(Q) \subset \mathbb{Z}_{+}^{m}$, so that
$\L(Q)\subset\{X_{1}^{a_{1}}\cdots X_{m}^{a_{m}}Z\,|\,a_{i}\geq
0\}$. After changing the lattice of reference
$\gp(S(P))=\mathbb{Z}^{d+1}$. Thus, for all $x\in L(cP)$,
$x=a_{0}x_{0}+\cdots+a_{n}x_{n}$ for some integers $a_{i}$ such that
$\sum_{i=0}^{n}{a_{i}}=c$. Let $f:k[cP] \to k[Q]$ be a graded
homomorphism.

By Section 3, $\varphi_{x}=f(x)Z^{-1}$ gives rise to the integral
affine map $\alpha:cP \to \mathbb{R}_{+}^{l}$ which respects the
binomial relations in $\L(cP)$. By Lemma 2 there exists a vector $v
\in \mathbb{Z}_{+}^{l}$ and an affine integral map $\beta:P\to
\mathbb{R}_{+}^{l}$ such that $\alpha=v+c\beta$. Now, $\beta(x_{i})$
gives rise to a polynomial $\theta_{i}\in k[X_{1},...\,,X_{m}]$ for
each $i\in\{0,1,\ldots n\}$. Similarly $v$ gives rise to a
polynomial $\psi\in k[X_{1},...\,,X_{m}]$. Then for $x\in L(cP)$,
$\varphi_{x}=\psi \theta_{0}^{a_{0}} \cdots \theta_{n}^{a_{n}}$.
Consider the $c^{th}$ homothetic blow-up of
\[
\Theta : k[P] \to k[Q']\,\,\mbox{such that}\,\,
\Theta(x_{i})=\theta_{i}Z
\]
for all $x_{i}\in L(P)$. $f$ is obtained by a polytope change
applied to $\Psi \star \Theta^{c}$ where
\[
\Psi : k[cP] \to k[Q'] \,\,\mbox{such that}\,\,\Psi(x)=\psi Z
\]
for all $x\in L(cP)$  and where $Q'$ is a large enough lattice
polytope to contain all relevant polytopes.

$\Theta$ is tame by assumption since it is a homomorphism from
$k[P]$. $\Psi$ factors as
\[
k[cP]\xrightarrow{\text{$x\mapsto
t$}}k[t]\xrightarrow{\text{$t\mapsto \psi Z$}}k[Q']
\]
where the first factor is tame since it maps monomials to monomials,
and the second factor is tame because it is a free extension of the
identity map $k\to k$. As a result, $f$ is tame.
\end{proof}

\section{Segre Products}\label{S:products}
The first step in proving Theorem $1(c)$ is Lemma $3$. Note that
$\L(P\times Q)=\L(P)\times L(Q)$. Suppose
$\L(P)=\{x_{0},\ldots,x_{n}\}$ and $\L(Q)=\{y_{0},\ldots,y_{m}\}$.
The copy of $P$ in $P\times Q$ with the lattice points
$\{(x_{0},y_{j}),\ldots,(x_{n},y_{j})\}$ is denoted by $P\times
y_{j}$. Similarly $x_{i}\times Q$ denotes the copy of $Q$ with
lattice points $\{(x_{i},y_{0}),\ldots,(x_{i},y_{m})\}$.

\begin{lemma}
Suppose $P$ and $Q$ are lattice polytopes and $\alpha: \L(P\times Q)
\to \mathbb{Z}^{l}_{+}$ is an affine map. Then $\alpha=a_{P}+a_{Q}$
for affine maps $a_{P},\, a_{Q}: \L(P\times Q) \to
\mathbb{Z}^{l}_{+}$ satisfying $a_{P}(P\times y_{i})=a_{P}(P\times
y_{j})$ and $a_{Q}(x_{i}\times Q)=a_{Q}(x_{j}\times Q)$ for all
$i,j$
\end{lemma}

\begin{proof}[Proof of Lemma 3]
Suppose $P$ and $Q$ are lattice polytopes and let $\alpha:
\L(P\times Q) \to \mathbb{Z}^{l}_{+}$ be an affine map. Since it is
enough to prove the lemma for each component of the maps involved,
there is no loss of generality to assume that $l=1$. Suppose
$\alpha$ takes its minimum value at the vertex $(x',y')\in\L(P\times
Q)$. Let $(x,y)\in\L(P\times Q)$. Sine $\alpha$ is affine,
\[
\alpha(x,y)+\alpha(x',y')=\alpha(x,y')+\alpha(x',y).
\]
In other words,
\[
\alpha(x,y)=\alpha(x,y')+\alpha(x',y)-\alpha(x',y').
\]
Define $a_{P}(x,y)=a(x,y')$ and $a_{Q}(x,y)=a(x',y)-a(x',y')$. By
the choice of $(x',y')$, the nonnegativity requirements are
satisfied.
\end{proof}
Now we can prove the desired result.

\begin{proof}[Proof of Theorem 1(c)]
Let $f: k[P~\times~Q]\to k[R]$ be a graded homomorphism, where $R$
is some lattice polytope. Assume $\L(P\times Q)=\L(P)\times
\L(Q)\subset \mathbb{R}_{+}^{d+e}$ such that $\L(P)\cap
\L(Q)=\{0\}$, and $\L(R)\subset \mathbb{R}_{+}^{c}$ (by applying a
polytope change if necessary). Suppose $\L(P)=\{x_{0},...\,
,x_{n}\}$ and $\L(Q)=\{y_{0},...\, ,y_{m}\}$.

Since $k[R]\subset \{X_{1}^{a_{1}}\cdot \cdot \cdot X_{c}^{a_{c}}Z |
\, a_{i}\geq 0\}$, the polynomials
$\varphi(x_{i},y_{j})=f(x_{i},y_{j})Z^{-1}$ belong to $k[X_{1},...\,
,X_{c}]$. Thus, by Section 3, $\varphi$ gives rise to an affine map
$\alpha:\L(P\times Q)\to \mathbb{Z}_{+}^{l}$ for some $l\in
\mathbb{N}$. Hence, by Lemma 3, $\alpha=a_{P}+a_{Q}$ where
$a_{P},a_{Q}: \L(P\times Q)\to\mathbb{Z}^{l}_{+}$ are affine maps
such that $a_{P}(x_{i},y_{j})=p_{i}$ and $a_{Q}(x_{i},y_{j})=q_{j}$
for some $p_{i},\,q_{j}\in\mathbb{Z}^{l}_{+}$.

Now, each $p_{i}$ and $q_{j}$ gives rise to polynomials
$\pi_{i},\,\varrho_{j}\in k[X_{1},...\, ,X_{c}]$, respectively. Then
$\varphi(x_{i},y_{j})=\pi_{i}\cdot\varrho_{j}$ for all
$(x_{i},y_{j})\in L(P\times Q)$. This means $f=f_{P}\star f_{Q}$
where
\[
\,\,\,f_{P}: k[P\times Q]\to k[R']\,\,\,\mbox{such that}\,\,
f_{P}(x_{i},y_{j})=\pi_{i}Z
\]
\[
\,\,\,f_{Q}: k[P\times Q]\to k[R']\,\,\,\mbox{such that}\,\,
f_{Q}(x_{i},y_{j})=\varrho_{j}Z
\]
for all $(x_{i},y_{j})\in \L(P\times Q)$ and where $R'$ is a large
enough lattice polytope to contain all relevant polytopes. $f_{P}$
is a tame homomorphism since it factors as
\[
k[P\times Q]\xrightarrow{\text{$(x_{i},y_{j})\mapsto
x_{i}$}}k[P]\xrightarrow{\text{$x_{i}\mapsto \pi_{i}Z$}}k[R']
\]
where the first factor is tame since it maps monomials to monomials
and the second is tame since it is a homomorphism from $k[P]$.
Similarly $f_{Q}$ is tame since it factors as
\[
k[P\times Q]\xrightarrow{\text{$(x_{i},y_{j})\mapsto
y_{j}$}}k[Q]\xrightarrow{\text{$y_{j}\mapsto \varrho_{j}Z$}}k[R']
\]
Therefore $f$ is a tame homomorphism, as desired.
\end{proof}

%

\bibliographystyle{amsalpha}
\bibliography{refTame}

\end{document}